\RequirePackage{fix-cm}
\documentclass[twocolumn]{svjour3}          
\smartqed  

\newcommand{\eg}{{\em e.g.}}
\newcommand{\ie}{{\em i.e.}}

\usepackage{graphicx}
\usepackage{rotating}
\usepackage{multirow}
\usepackage{multicol}
\usepackage{siunitx}
\usepackage{natbib}
\usepackage{epsfig} 
\usepackage{amsfonts}
\usepackage{amssymb,amsmath}
\usepackage{bm}
\usepackage[lofdepth,lotdepth]{subfig}
\usepackage{float}
\usepackage{multirow}
\usepackage{makecell}
\usepackage{algorithm,algpseudocode}
\algnewcommand{\Inputs}[1]{%
  \State \textbf{Inputs:}
  \Statex \hspace*{\algorithmicindent}\parbox[t]{.8\linewidth}{\raggedright #1}
}
\algnewcommand{\Initialize}[1]{%
  \State \textbf{Initialize:}
  \Statex \hspace*{\algorithmicindent}\parbox[t]{.8\linewidth}{\raggedright #1}
}

\usepackage{algorithm}
\usepackage{algpseudocode}
\algnewcommand{\LineComment}[1]{\State \(\triangleright\) #1}

\DeclareMathOperator*{\argmax}{arg\,max}

\begin{document}
\sloppy

\title{Adaptive Expansion Bayesian Optimization for Unbounded Global Optimization\thanks{
This work was funded by The Defense Advanced Research Projects Agency (DARPA-16-63-YFA-FP-059) via the Young Faculty Award (YFA) Program.
}
}


\author{Wei Chen         \and
        Mark Fuge 
}


\institute{W. Chen \at
              Department of Mechanical Engineering \\
              University of Maryland \\
              College Park, MD 20742 \\
              \email{wchen459@umd.edu}           
           \and
           M. Fuge \at
              Department of Mechanical Engineering \\
              University of Maryland \\
              College Park, MD 20742
}

\date{Received: date / Accepted: date}

\maketitle

\begin{abstract}
  Bayesian optimization is normally performed within fixed variable bounds. In cases like hyperparameter tuning for machine learning algorithms, setting the variable bounds is not trivial. It is hard to guarantee that any fixed bounds will include the true global optimum. We propose a Bayesian optimization approach that only needs to specify an initial search space that does not necessarily include the global optimum, and expands the search space when necessary. However, over-exploration may occur during the search space expansion. Our method can adaptively balance exploration and exploitation in an expanding space. Results on a range of synthetic test functions and an MLP hyperparameter optimization task show that the proposed method out-performs or at least as good as the current state-of-the-art methods.
\end{abstract}

\section{Introduction}

Bayesian optimization (BO) is a global optimization technique targeted for expensive black-box functions~\citep{shahriari2015taking}. Particularly, one of its important applications in the machine learning community is automated hyperparameter tuning~\citep{snoek2012practical,swersky2013multi,springenberg2016bayesian}. In a standard BO process, the objective function is modeled as a random function with a prior distribution. This prior updates to form a posterior after new observations (\ie, a Gaussian process or GP~\citep{rasmussen2006gaussian}). The decision about which observation to collect next is made by globally maximizing an acquisition function based on the posterior. This step requires fixed variable bounds, which are sometimes not trivial to set. It is hard to guarantee that any fixed bounds will include the true global optimum.

In this paper, we modified the standard BO approach so that the fixed variable bounds are not required. When the search space is unbounded, the acquisition function can have suprema at infinity, where the uncertainty is maximized. Thus we search only in the region with sufficiently low uncertainty, which we referred to as the \textit{low-uncertainty region}. This low-uncertainty region expands as we add more observations. We call this method Adaptive Expansion Bayesian Optimization (AEBO). 

The main technical contributions of this paper are:
\begin{enumerate}
\item An acquisition strategy that bounds the GP model uncertainty and  adaptively expands the low-uncertainty region; and
\item Theoretical results regarding how to adaptively set the threshold of the uncertainty bound to avoid the over-exploration problem that occurs in an expanding search space.
\end{enumerate}

\section{Bayesian Optimization}

Bayesian Optimization uses a sequential strategy to search for the global optimum of expensive black-box functions. Assuming we have an objective function: $f: \mathbb{R}^d\rightarrow\mathbb{R}$, and the observation of its output has Gaussian noise: $y\sim\mathcal{N}(f(\mathbf{x}), \sigma_n^2)$. It is expensive to evaluate either the function $f$ or its gradient (assuming we can only approximate the gradient by the finite difference method when $f$ is a black-box function and that Automatic Differentiation methods cannot be used). Thus the goal of BO is to minimize the number of evaluations needed to find the global minimum solution. BO treats the objective as a random function that has a prior distribution, and update this prior to form a posterior distribution over the function after observing data. This can be done by using a Gaussian process (GP). The posterior distribution can then be used to form an acquisition criterion that proposes to evaluate $f$ at a promising point so that the regret is minimized. The GP posterior can then be updated after the new observation. This process repeats until the evaluation budget runs out or a satisfied solution is achieved.
We will elaborate on the Gaussian process and the acquisition function in the following sections.

\subsection{Gaussian Process}

The Gaussian process (GP)~\citep{rasmussen2006gaussian} estimates the distribution of the objective function. A kernel (covariance) function $k(\mathbf{x},\mathbf{x}')$ is used to measure the similarity between two points $\mathbf{x}$ and $\mathbf{x}'$. It encodes the assumption that ``similar inputs should have similar outputs". The specific choice of the kernel is not central to the core contributions of the paper. 

Given $N$ observations $\mathcal{D} = (\mathbf{X}, \mathbf{y}) = \{(\mathbf{x}_i, y_i)|i=1,...,N\}$, the GP posterior $f(\mathbf{x})$ at any point $\mathbf{x}$ is a Gaussian distribution: $f(\mathbf{x})|\mathcal{D},\mathbf{x} \sim \mathcal{N}(\mu(\mathbf{x}), \sigma^2(\mathbf{x}))$ with the mean and the variance expressed as
\begin{align}
\mu(\mathbf{x}) &= \mathbf{k}(\mathbf{x})^\top(\mathbf{K}+\sigma_n^2\mathbf{I})^{-1}{\mathbf{y}}
\label{eq:mean}\\
\sigma^2(\mathbf{x}) &= k(\mathbf{x},\mathbf{x})-\mathbf{k}(\mathbf{x})^\top(\mathbf{K}+\sigma_n^2\mathbf{I})^{-1}\mathbf{k}(\mathbf{x})
\label{eq:variance}
\end{align}
where $\mathbf{k}(\mathbf{x})$ is an $N$-dimensional vector with the $i$-th dimension being $k(\mathbf{x},\mathbf{x}_i)$, and $\mathbf{K}$ is an $N\times N$ covariance matrix with $K_{ij}=k(\mathbf{x}_i,\mathbf{x}_j)$.

\subsection{Acquisition Function}

Bayesian optimization picks the next point to evaluate by maximizing an acquisition function, which is computed based on the GP posterior. Common acquisition functions include the probability of improvement (PI)~\cite{kushner1964new}, the expected improvement (EI)~\cite{jones1998efficient}, the Gaussian Process upper confidence bound (GP-UCB)~\cite{srinivas2009gaussian}, and those based on entropy search~\cite{hennig2012entropy,hernandez2014predictive,wang2017max}. 

In this paper, we use EI as our acquisition function. It measures the expected amount of improvement over the current best solution based on the learned GP model:
\begin{equation}
\begin{aligned}
\text{EI}(\mathbf{x}) &= \mathbb{E}[\max\{0,f(\mathbf{x})-f'\}] \\
&= \int_{f'}^{+\infty}(f-f')\mathcal{N}(f;\mu(\mathbf{x}), \sigma^2(\mathbf{x}))df \\
&= \sigma(\mathbf{x})(u\Phi(u) + \phi(u))
\end{aligned}
\label{eq:ei}
\end{equation}
where $f'$ is the current best objective function value, $u=(\mu(\mathbf{x})-f')/\sigma(\mathbf{x})$, and $\Phi$ and $\phi$ are the cumulative density function (CDF) and probability density function (PDF) of the standard normal distribution, respectively.

\subsection{Previous Work on Unbounded Bayesian Optimization}

Normally Bayesian optimization is performed within fixed variable bounds. But in cases such as algorithm hyperparameter tuning~\citep{shahriari2016unbounded,swersky2013multi,springenberg2016bayesian} and shape optimization~\citep{palar2019efficient}, setting the variable bounds are not trivial. It is hard to guarantee that any fixed bounds will include the true global optimum.

Two types of solutions were proposed to handle this problem: 1)~performing BO in an unbounded space by regularization via non-stationary prior means so that the acquisition function's suprema will not be at infinity (\textit{regularization-based methods})~\citep{shahriari2016unbounded,rainforth2016bayesian}; and 2)~performing BO in ``soft bounds'' that are gradually expanded over iterations (\textit{expansion-based methods})~\citep{shahriari2016unbounded,nguyen2017bayesian,nguyen2018filtering}. The first solution computes an acquisition function that is biased toward regions near some user-specified center point, thus insufficient exploitation may occur when the optimal solution is far from the center. The second solution either expands each direction equally, which often yields to unnecessarily large search spaces~\citep{shahriari2016unbounded,rainforth2016bayesian}; or expands only to the promising region where the upper confidence bound (UCB) is larger than the lower confidence bound (LCB) of the current best solution~\citep{nguyen2017bayesian,nguyen2018filtering}. But the latter approach has to perform global search twice at each iteration\textemdash one for the maximum LCB to filter out the non-promising region; and one for the maximum acquisition function value. It expands the bounds of the first search according to a hard-coded rule, and hence may show a lack of adaptability to different optimization problems.

The adaptive expansion Bayesian optimization is different from the previous methods in that it adaptively expands the search space based on the uncertainty of the GP model. It can essentially avoid the aforementioned issues by employing a strategy that we will introduce in the following section.

\section{Adaptive Expansion Bayesian Optimization}

In this section, we will introduce the main ingredients of AEBO, namely, its acquisition strategy (Sect.~\ref{sec:acquisition_strategy}), global optimization of the acquisition function (Sect.~\ref{sec:feasible_domain_bounds}), the way of adaptively balancing exploration and exploitation (Sect.~\ref{sec:explore_exploit}), and a trick to improve exploitation when expanding the search space (Sect.~\ref{sec:local_search}).

\subsection{Acquisition Strategy}
\label{sec:acquisition_strategy}

Our acquisition strategy can be expressed as the following constrained optimization problem:
\begin{equation}
\begin{aligned}
\max_{\mathbf{x}\in\mathbb{R}^d} ~& \text{EI}(\mathbf{x}) \\
\text{s.t.} ~& \sigma^2(\mathbf{x}) \leq \tau k_0
\end{aligned}
\label{eq:strategy}
\end{equation}
where $\tau\in (0,1)$ is a coefficient controlling the aggressiveness/conservativeness of exploration, and $k_0=\sigma^2(\mathbf{x}_\infty)$ with $\mathbf{x}_\infty$ denotes a point infinitely far away from the observations. Based on Eq.~\ref{eq:variance}, we have $k_0=\sigma^2(\mathbf{x}_\infty)=k(\mathbf{x},\mathbf{x})$. When using a RBF or a Mat\'ern kernel, for example, simply we have $k_0=1$. Under this acquisition strategy, only points with low GP model uncertainty will be picked for evaluation. When a GP adds evaluations, the $\sigma^2(\mathbf{x}) \leq \tau k_0$ region will always expand, because the predictive variance near the added evaluations will decrease. Thus our method guarantees a growing search space.

To avoid excessive local exploitation, we can modify the equation for EI (Eq.~\ref{eq:ei}) as
\begin{equation}
\text{EI}(\mathbf{x}) = \mathbb{E}[\max\{0,f(\mathbf{x})-(f'+\epsilon)\}]
\label{eq:epsilon}
\end{equation}
where $\epsilon>0$ is the \textit{minimum improvement parameter}~\citep{jones2001taxonomy,shahriari2016unbounded,shahriari2016practical}.

In many real-world cases, we need to deal with constrained Bayesian optimization problems of the following two kinds: 1)~there are infeasible regions in the input space (\eg, some experimental configurations are infeasible); and 2)~the objective function $f$ does not have definition in some regions of the input space (\eg, when the hyperparameters of a neural network are not properly chosen, exploding gradients may occur, which may lead to NaN weight values and hence the NaN accuracy). This is especially common when the we have an unbounded or expanding search space. Therefore, it is worth extending AEBO to make it suitable for constrained BO problems. Specifically, we can modify Eq.~\ref{eq:strategy} based on Refs.~\citep{basudhar2012constrained} and \citep{gelbart2014bayesian}:
\begin{equation}
\begin{aligned}
\max_{\mathbf{x}\in\mathbb{R}^d} ~& \text{EI}(\mathbf{x})\text{Pr}(\mathcal{C}(\mathbf{x})) \\
\text{s.t.} ~& \sigma^2(\mathbf{x}) \leq \tau k_0 \\
& \text{Pr}(\mathcal{C}(\mathbf{x})) \geq 0.5
\end{aligned}
\label{eq:constrained_strategy}
\end{equation}
where $\mathcal{C}(\mathbf{x})$ is an indicator of whether the constraints are satisfied or whether the objective function has definition.

\subsection{Feasible Domain Bounds}
\label{sec:feasible_domain_bounds}

\begin{figure}[h]
\centering
\includegraphics[width=0.25\textwidth]{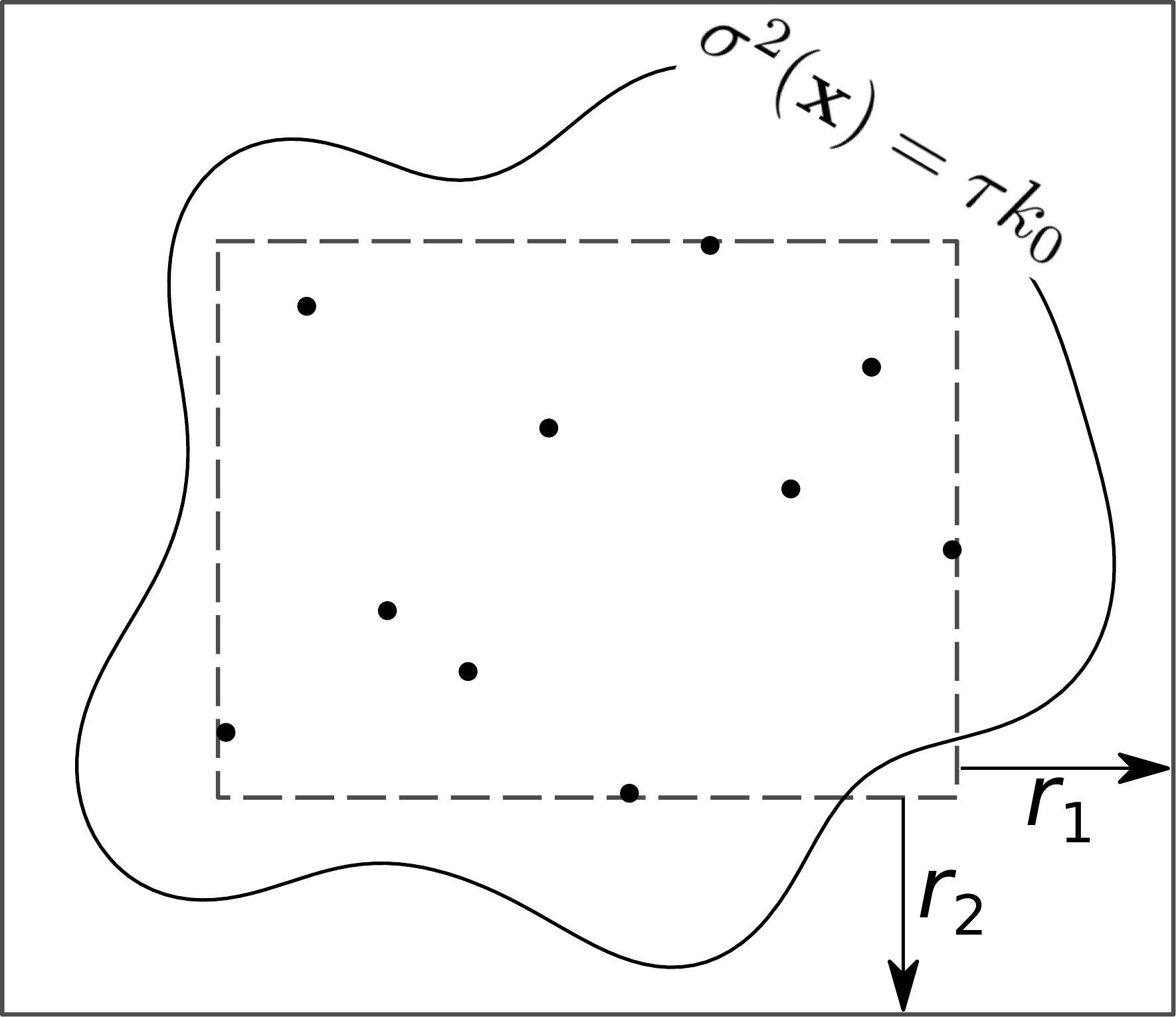}
\caption{Feasible domain bounds. In each iteration, we expand the minimum bounding box of evaluated samples along the $i$-th axis by $r_i$.}
\label{fig:bounds}
\end{figure}


The feasible domain defined by Eq.~\ref{eq:strategy} is bounded by the isocontour $\sigma^2(\mathbf{x})=\tau k_0$. However, it is easier to search inside a bounding box instead of an irregular isocontour when solving the global optimization problem in Eq.~\ref{eq:strategy}. We can show that the solution to Eq.~\ref{eq:strategy} is inside a bounding box, which we call the \textit{feasible domain bounds}. The feasible domain bounds can be derived by expanding the minimum bounding box of evaluated samples along the $i$-th axis by an expansion rate $r_i$ (Fig.~\ref{fig:bounds}).
Then constrained global optimization of the acquisition function can be performed within that feasible domain bounds. The derivation of $r_i$ is included in the supplementary material.

\subsection{Adaptive Exploration-Exploitation Trade-off}
\label{sec:explore_exploit}

\begin{figure*}
\centering
\includegraphics[width=1\textwidth]{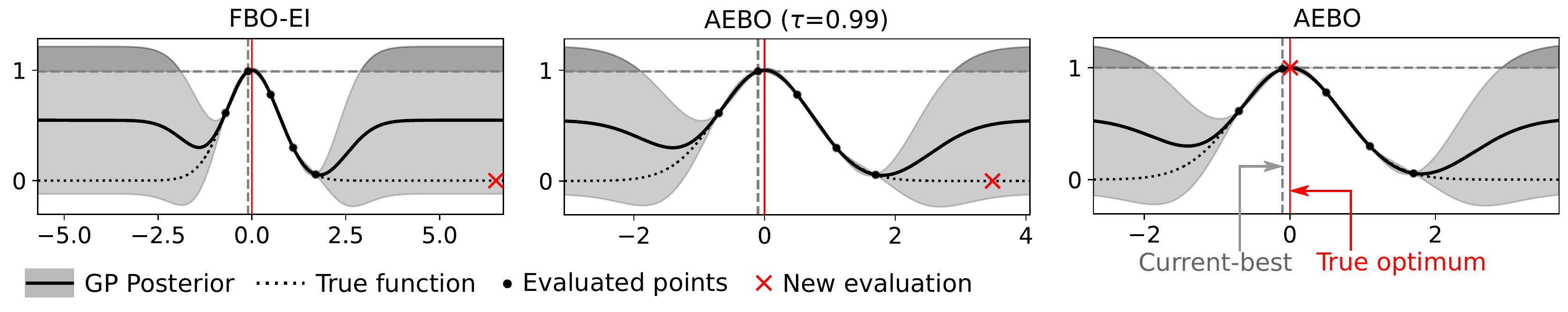}
\caption{Over-exploration in FBO-EI (left) and AEBO with large $\tau$ (middle); By adaptively setting $\tau$, AEBO enforces exploitation based on an accuracy criterion (right).}
\label{fig:over_explore_1d}
\end{figure*}

\begin{figure*}
\centering
\includegraphics[width=0.75\textwidth]{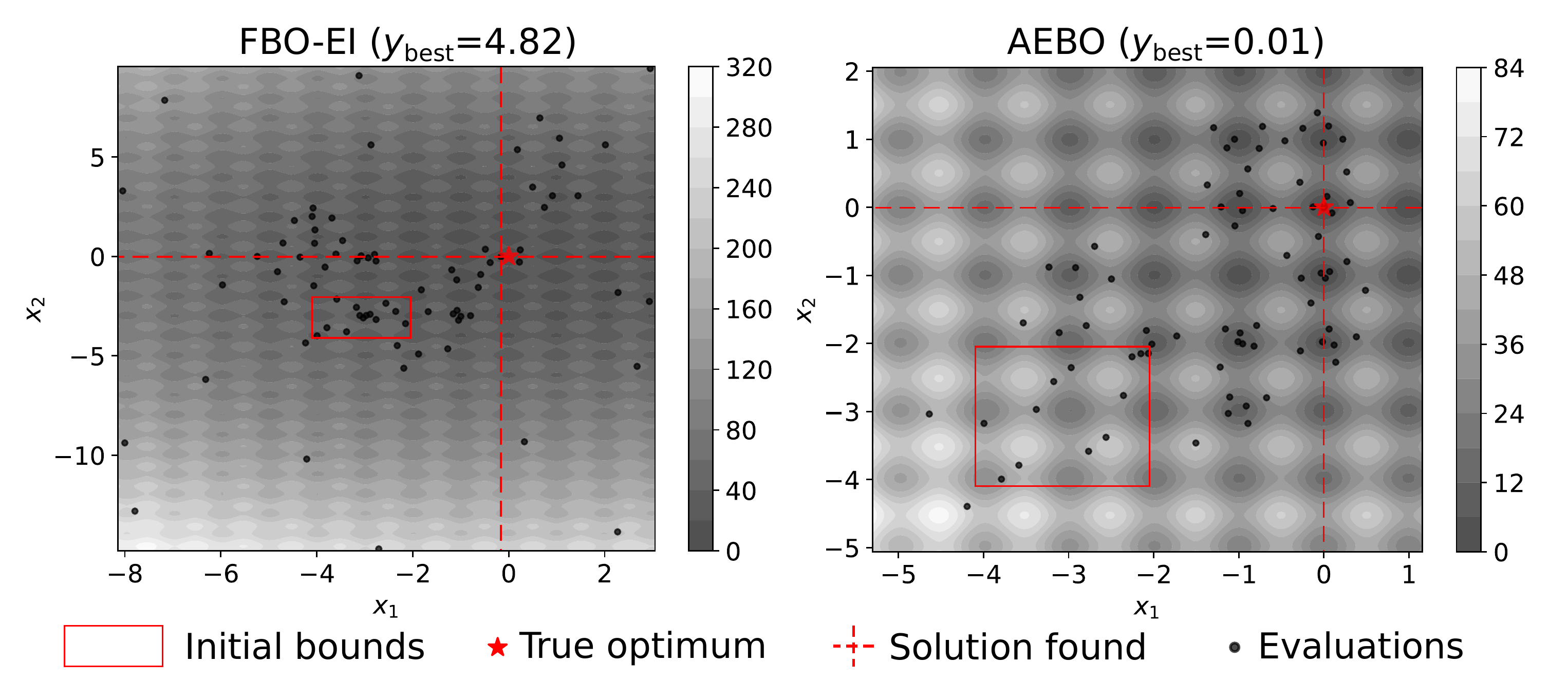}
\caption{In contrast with AEBO (right), FBO-EI (left) spends too much budget on randomly exploring the search space (left).}
\label{fig:over_explore_2d}
\end{figure*}

A problem of an expanding search space is that a new evaluation may get too far away from the region of interest, due to the high uncertainty and hence the high EI in far-away regions (Fig.~\ref{fig:over_explore_1d}). The informativeness of those high-uncertainty regions, however, is low because of the sparsity of observed data near them. Thus sampling at those regions is like shooting in the dark. This is fine in the fixed-bound BO, because it will exploit the region of interest eventually after finishing exploring those high-uncertainty regions. However, with expanding bounds, the high-uncertainty regions are expanding and BO could continuously sample in those regions and never head back to exploit the region of interest. As a result, the algorithm will spend too much budget on randomly exploring the search space but have insufficient exploitation, as shown in the left plot of Fig.~\ref{fig:over_explore_2d}. We call this \textit{over-exploration}. 

This over-exploration problem can exist in every unbounded Bayesian optimization algorithm with an aggressive expansion strategy. Adaptive Expansion BO can solve this problem by avoiding exploring in regions where our estimated model is uncertain (\ie, constraining the GP's predictive variance $\sigma^2(\mathbf{x})$, see Eq.~\ref{eq:strategy}). However, one has to choose a proper coefficient $\tau$ to set the uncertainty threshold. In this section, we derive a way of setting $\tau$ adaptively to balance exploration and exploitation as the search space expands.

The simplest way to avoid over-exploration is to force the algorithm to stop exploring and start to refine the solution by exploiting near the current best point. In AEBO, exploration is performed by sampling points along the feasible domain boundary (\ie, $\sigma^2(\mathbf{x})=\tau k_0$). Thus, we can avoid over-exploration by decreasing $\tau$ so that the expected improvement on the boundary is lower than that near the current best solution $\mathbf{x}'$.

The expected improvement on the boundary can be expressed via the predictive mean $\mu_{\tau}$:
\begin{equation}
\text{EI}_{\tau}(\mu_{\tau}) = (\mu_{\tau}-f')\Phi\left(\frac{\mu_{\tau}-f'}{\sqrt{\tau k_0}}\right) + \sqrt{\tau k_0}\phi\left(\frac{\mu_{\tau}-f'}{\sqrt{\tau k_0}}\right)
\label{eq:ei_tau}
\end{equation}
Also we have
\begin{equation}
\max_{\sigma^2(\mathbf{x})=\tau k_0}\{\text{EI}_{\tau}(\mu_\tau)\} = \text{EI}_{\tau}\left(\max_{\sigma^2(\mathbf{x})=\tau k_0}\{\mu_\tau\}\right) = \text{EI}_\tau(\mu_m)
\label{eq:max_ei_tau}
\end{equation}
since EI monotonically increases with the predictive mean.

The expected improvement near the current best solution $\mathbf{x}'$ is
\begin{equation}
\text{EI}_+ = (\mu_+-f')\Phi\left(\frac{\mu_+-f'}{\sigma_+}\right) + \sigma_+\phi\left(\frac{\mu_+-f'}{\sigma_+}\right)
\label{eq:ei_plus}
\end{equation}
where $\mu_+$ and $\sigma_+$ are the predictive mean and standard deviation respectively at a point $\mathbf{x}_+$ near the current best solution. Assuming that the GP mean function $\mu(\mathbf{x})$ is Lipschitz continuous, we have $f'-\mu_+=\delta$, where $\delta$ is a small positive real number. Thus we have $u_+=-\delta/\sigma_+$.

Now we can set $\text{EI}_+ > \text{EI}_{\tau}(\mu_m)$ to encourage exploitation. However, we do not want pure exploitation. Specifically, we want to stop exploitation at $\mathbf{x}_+$ whenever the \textit{room for improvement} over the current solution $f'$ within the neighborhood of $\mathbf{x}'$ is sufficiently low with a high probability:
\begin{equation}
\text{Pr}(f_+-f'\leq\xi) \geq 1-\kappa
\label{eq:improve_prob}
\end{equation}
where $f_+ \sim \mathcal{N}(\mu_+, \sigma_+^2)$, and $\xi\geq 0$ and $0<\kappa<1$ are small real numbers. 
From Eq.~\ref{eq:improve_prob} we can derive
\begin{equation*}
\sigma_+ \leq \frac{\xi+\delta}{\Phi^{-1}(1-\kappa)} = \sigma_0
\label{eq:sigma_plus}
\end{equation*}
Thus we only need to exploit at $\mathbf{x}_+$ when $\sigma_+ > \sigma_0$. By substituting it into Eq.~\ref{eq:ei_plus}, we get a lower bound for $\text{EI}_+$:
\begin{equation}
\text{EI}_+ > -\delta\Phi(-\delta/\sigma_0) + \sigma_0\phi(-\delta/\sigma_0) = \text{EI}_0
\label{eq:ei_plus_lb}
\end{equation}
since EI monotonically increases with the predictive variance. We can set this lower bound $\text{EI}_0$ equal to $\text{EI}_{\tau}(\mu_m)$ to enforce $\text{EI}_+ > \text{EI}_{\tau}(\mu_m)$ when $\sigma_+ > \sigma_0$ (\ie, when exploitation is necessary). Using Eq.~\ref{eq:ei_tau}, we can write $\text{EI}_{\tau}(\mu_m)=\text{EI}_0$ as
\begin{equation}
(\mu_m-f')\Phi\left(\frac{\mu_m-f'}{\sqrt{\tau k_0}}\right) + \sqrt{\tau k_0}\phi\left(\frac{\mu_m-f'}{\sqrt{\tau k_0}}\right) = \text{EI}_0
\label{eq:set_tau}
\end{equation}
We can solve for $\tau$ by using any root finding algorithm (\eg, Newton's method).

At the beginning of the optimization process, we do not need to make sure the room for improvement over $f'$ is small within the neighborhood of $\mathbf{x}'$. Rather, we want to explore other regions that may contain better local optima. Thus we can set $\xi=\xi_0$ at the beginning, where a larger $\xi_0$ allows more exploration, and then linearly anneal $\xi$ over iterations until $\xi=0$ (\eg, towards the end of a computational budget). As a result, AEBO's focus gradually switches from exploration to exploitation.

In practice, we can set $\mu_m$ as the prior mean (0 by default), since it is usually the case when over-exploration occurs. Thus this adaptive approach can effectively avoid over-exploration, as shown in Fig.~\ref{fig:over_explore_1d} and Fig.~\ref{fig:over_explore_2d}. If in reality $\mu_m<0$, then $\text{EI}_{\tau}(\mu_m) < \text{EI}_0 < \text{EI}_+$, AEBO will exploit near the current best solution even when it is unnecessary (\ie, $\sigma_+ < \sigma_0$); while if $\mu_m>0$, then $\text{EI}_{\tau}(\mu_m) > \text{EI}_0$, AEBO may explore when it should exploit.

\subsection{Local Search for Better Exploitation}
\label{sec:local_search}

In practice, to perform the global optimization of Eq.~\ref{eq:strategy}, we can sample initial \textit{candidate solutions} within the bounds derived in Sect.~\ref{sec:feasible_domain_bounds}, and refine those solutions using a constrained optimization method (\eg, Sequential Least Squares Programming). If the GP kernel is fixed, the search space is always expanding, because $\sigma$ is monotonically non-increasing as the number of observations increases. Specifically, the space near the queried point will be added to the search space volume. This results in a volume increase that is exponential with respect to the search space dimensionality. It will become harder for the candidate solutions to maintain the coverage of the search space as the optimization proceeds, especially when the problem has high dimensionality. Although we keep increasing exploitation by annealing $\xi$, it does not guarantee that we will exploit near the current best solution in a large search space. This problem was not addressed in previous unbounded Bayesian optimization methods~\citep{shahriari2016unbounded,nguyen2018filtering}. A straight-forward way to solve the problem is to increase the density of search algorithms, but this continuously increases the computational cost for each iteration. Alternatively, we propose local search near the current best solution to allow better exploitation. Specifically, in each iteration, we generate the same number of candidate solutions but divide it for two tasks\textemdash global search and local search. Global search tries to find a promising point in the entire feasible domain in Eq.~\ref{eq:strategy}; while local search tries to find a promising point near the current best solution. This avoids insufficient exploitation but will not increase the computational cost.

The optimization process is summarized in Algorithm~\ref{alg:aebo}.

\begin{algorithm}
\caption{Adaptive Expansion Bayesian optimization}
\label{alg:aebo}
\begin{algorithmic}[1]
    \LineComment{Given objective function $f$, initial bounds $\mathcal{B}$, initial evaluation $n$, and evaluation budget $N$}
\Procedure{Maximize}{$f, \mathcal{B}, n, N$}
    \State Sample $n$ points $\{\mathbf{x}_1, ..., \mathbf{x}_n\}$ in $\mathcal{B}$ using LHS
    \State $y_i \leftarrow f(\mathbf{x}_i)$, $\forall i=1,...,n$
    \State $\mathcal{D} \leftarrow \{(\mathbf{x}_1, y_1), ..., (\mathbf{x}_n, y_n)\}$
    \State $f' \leftarrow \max_i\{y_i\}$, $\mathbf{x}^* \leftarrow \argmax_{\mathbf{x}_i}\{y_i\}$
    \For {$t = (n+1):N$}
    	\State Fit the GP model $\mathcal{M}$ to $\mathcal{D}$
    	\State Compute $\tau$ based on Eq.~\ref{eq:set_tau}
    	\State Expand the minimum bounding box of $\mathcal{D}$ by $r_i$ to get the feasible domain bounds $\mathcal{B}'$
    	\State Search for the solution $\mathbf{x}_t$ to Eq.~\ref{eq:strategy} inside $\mathcal{B}'$
    	\State $y_t \leftarrow f(\mathbf{x}_t)$
    	\If{$y_t > f'$}
    	    \State $\mathbf{x}^* \leftarrow \mathbf{x}_t$, $f' \leftarrow y_t$
    	\EndIf
    	\State $\mathcal{D} \leftarrow \mathcal{D}\bigcup \{(\mathbf{x}_t, y_t)\}$
    \EndFor
	\State\Return $\mathbf{x}^*$
\EndProcedure
\end{algorithmic}
\end{algorithm}

\begin{table*}[h]
\centering
\caption{Optimization results for synthetic benchmarks}
\label{tab:test_func_results}
\begin{tabular}{lrrrrrrrr}
\hline\noalign{\smallskip}
Method & SixHumpCamel & Branin & Rastrigin & Hartmann3  \\
\noalign{\smallskip}\hline\noalign{\smallskip}
AEBO & $\mathbf{-1.03 \pm 0.00}$ & $\mathbf{0.40 \pm 0.00}$ & $\mathbf{0.26 \pm 0.43}$ & $\mathbf{-3.69 \pm 0.22}$ \\
FBO-EI & $-0.97 \pm 0.05$ & $1.27 \pm 1.49$ & $4.52 \pm 2.29$ & $-1.50 \pm 0.93$ \\
FBO-UCB & $-0.85 \pm 0.09$ & $1.43 \pm 0.61$ & $5.50 \pm 2.60$ & $-2.31 \pm 1.34$ \\
EI-Q & $-0.28 \pm 0.37$ & $2.95 \pm 1.73$ & $8.10 \pm 1.47$ & $-2.43 \pm 0.65$ \\
EI-H & $-0.47 \pm 0.46$ & $1.89 \pm 1.00$ & $7.39 \pm 1.29$ & $-3.41 \pm 0.25$ \\
\noalign{\smallskip}\hline
\hline\noalign{\smallskip}
Method & Hartmann6 & Beale & Rosenbrock  \\
\noalign{\smallskip}\hline\noalign{\smallskip}
AEBO & $-3.29 \pm 0.03$ & $\mathbf{0.18 \pm 0.26}$ & $\mathbf{0.68 \pm 0.78}$ \\
FBO-EI & $\mathbf{-3.30 \pm 0.03}$ & $0.41 \pm 0.32$ & $7.72 \pm 9.25$ \\
FBO-UCB & $-3.26 \pm 0.04$ & $0.46 \pm 0.37$ & $17.39 \pm 33.04$ \\
EI-Q & $-2.32 \pm 0.23$ & $4.25 \pm 2.83$ & $17.45 \pm 20.12$ \\
EI-H & $-2.82 \pm 0.15$ & $3.87 \pm 3.31$ & $20.63 \pm 19.71$ \\
\noalign{\smallskip}\hline
\end{tabular}
\end{table*}

\begin{figure*}[h]
\centering
\includegraphics[width=1.0\textwidth]{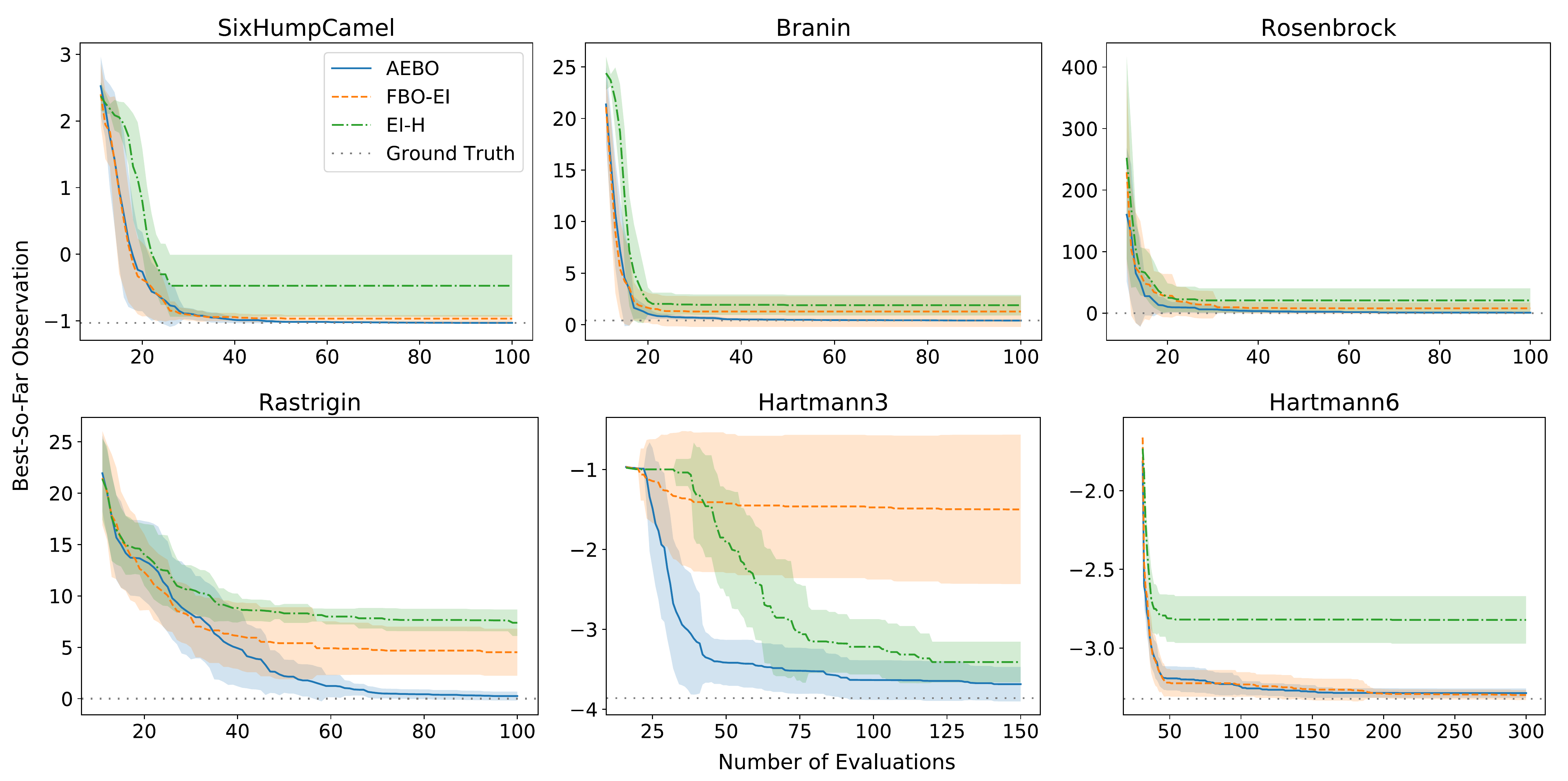}
\caption{Optimization history for synthetic benchmarks}
\label{fig:synthetic_history}
\end{figure*}

\section{Experiments}

We evaluate the AEBO on both a range of synthetic test functions and an MLP hyperparameter optimization task. We also demonstrate the effect of dimensionality on AEBO's performance and the significance of adaptive exploration-exploitation trade-off in AEBO. 

\subsection{Experimental Protocol}

The evaluation budget was set to $50d$, and an initial sample size of $5d$ was drawn by using Latin hypercube sampling~\citep{mckay1979comparison,jones1998efficient}. For simplicity we used an isotropic RBF kernel for the GP (\ie, $k(\mathbf{x},\mathbf{x}')=\exp\left(-\sum_{i=1}^d (\mathbf{x}_i-\mathbf{x}'_i)^2/(2l_i^2)\right)$, where $l_1=...=l_d$). We normalized the observed function outputs before fitting a GP regression model. We set $\xi_0=0.1$, $\kappa=0.1$, $\epsilon=0.01$, and $\delta = 0.01$. For each test function, we set the initial bounds to be [10\%, 30\%] of its original bounds, as was also configured in Ref.~\citep{nguyen2018filtering}. All the initial bounds do not include global optima. We compared AEBO to methods from Ref.~\citep{shahriari2016unbounded} (\ie, EI-Q and EI-H) and Ref.~\citep{nguyen2018filtering} (\ie, FBO-EI and FBO-UCB).

\subsection{Synthetic Benchmarks}

We used seven standard global optimization test functions. As shown in Table~\ref{tab:test_func_results}, AEBO out-performs other methods on most test functions. Note that AEBO's results have lower variance compared to other methods, which is an indication of robustness. Figure~\ref{fig:synthetic_history} shows the optimization history on benchmark functions. It shows that compared to the other two state-of-the-art methods, AEBO converged faster and achieved a better solution in most cases.

\begin{figure}[h]
\centering
\includegraphics[width=0.5\textwidth]{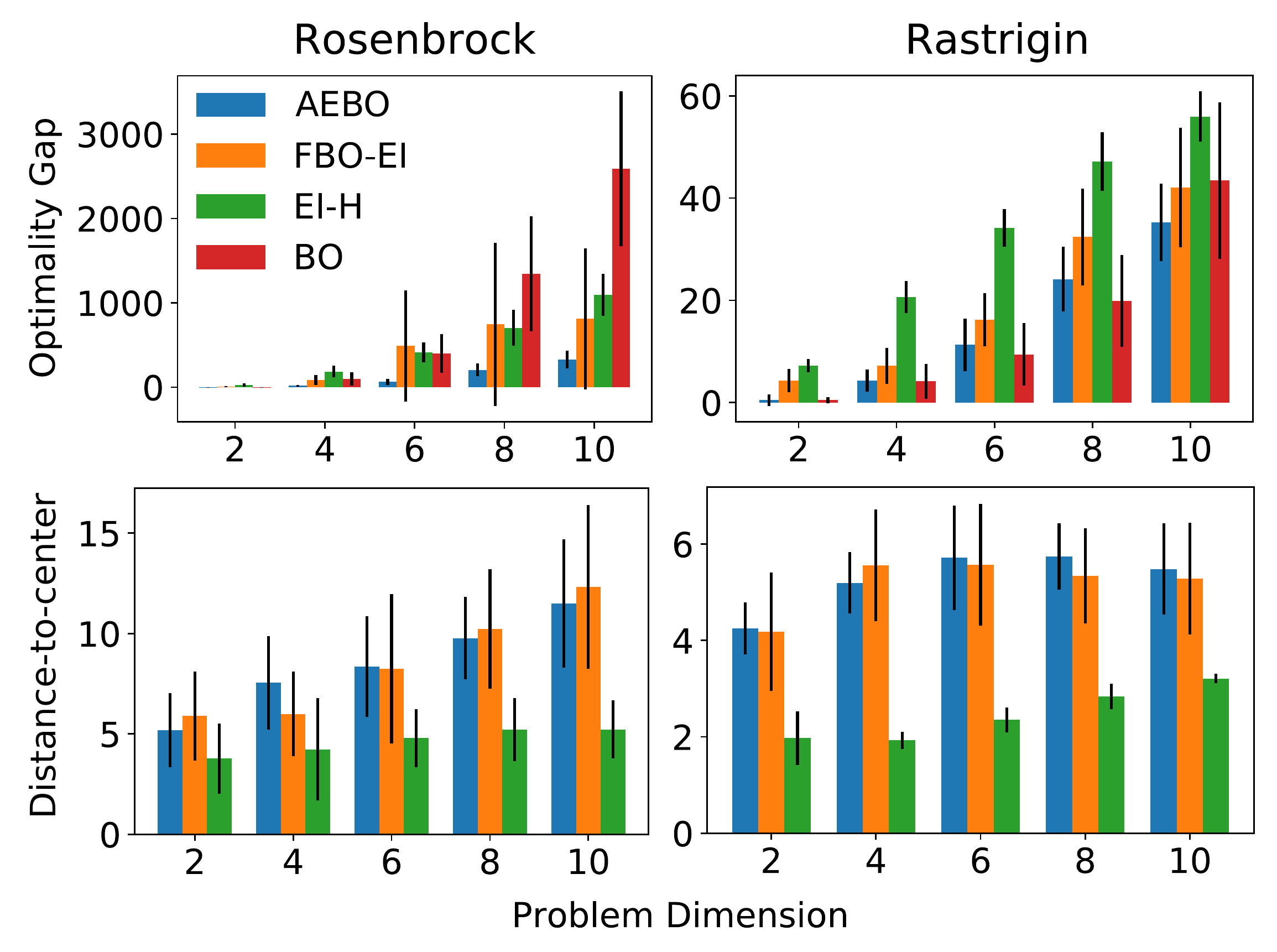}
\caption{The effect of problem dimension on the optimality gap (top) and the reachability of optimal solutions found by different methods (bottom). Note that the standard BO, unlike other methods, was performed within the bounds that include the actual global optima.}
\label{fig:dim_gap_dist}
\end{figure}

We demonstrated the effects of problem dimensionality on AEBO's performance by using two synthetic benchmarks, as shown in Fig.~\ref{fig:dim_gap_dist}. We compared AEBO to: 1)~the other two state-of-the-art methods\textemdash FBO-EI and EI-H, and 2)~the standard BO with the original function bounds (which include global optima). Here we define the \textit{optimality gap} $e = y_{sol}-y_{opt}$, where $y_{sol}$ and $y_{opt}$ are the minimal observation and the true minimum of the objective function, respectively. The results show that the optimality gap increases with the problem dimension, which can be explained by the curse of dimensionality~\citep{bellman1957dynamic}. AEBO demonstrated the best performance among methods dealing with an unbounded or expanding search space, and is almost as good as the standard BO on the Rastrigin function. Since the global optimum of the Rosenbrock function is inside a narrow flat valley, it is trivial to find the valley but difficult to converge to the global optimum. Thus it requires large budget for exploitation in that valley (\ie, \textit{exploitation-intense}). The standard BO and expansion-based methods like FBO-EI may have unnecessarily large search space and hence waste budget on exploring regions far from the global optimum, rather than exploiting the valley. Thus compared to AEBO, it is more difficult for these three methods to find good solutions on the Rosenbrock function, especially when the dimensionality is high (Fig.~\ref{fig:dim_gap_dist}).

\begin{figure}[h]
\centering
\includegraphics[width=0.5\textwidth]{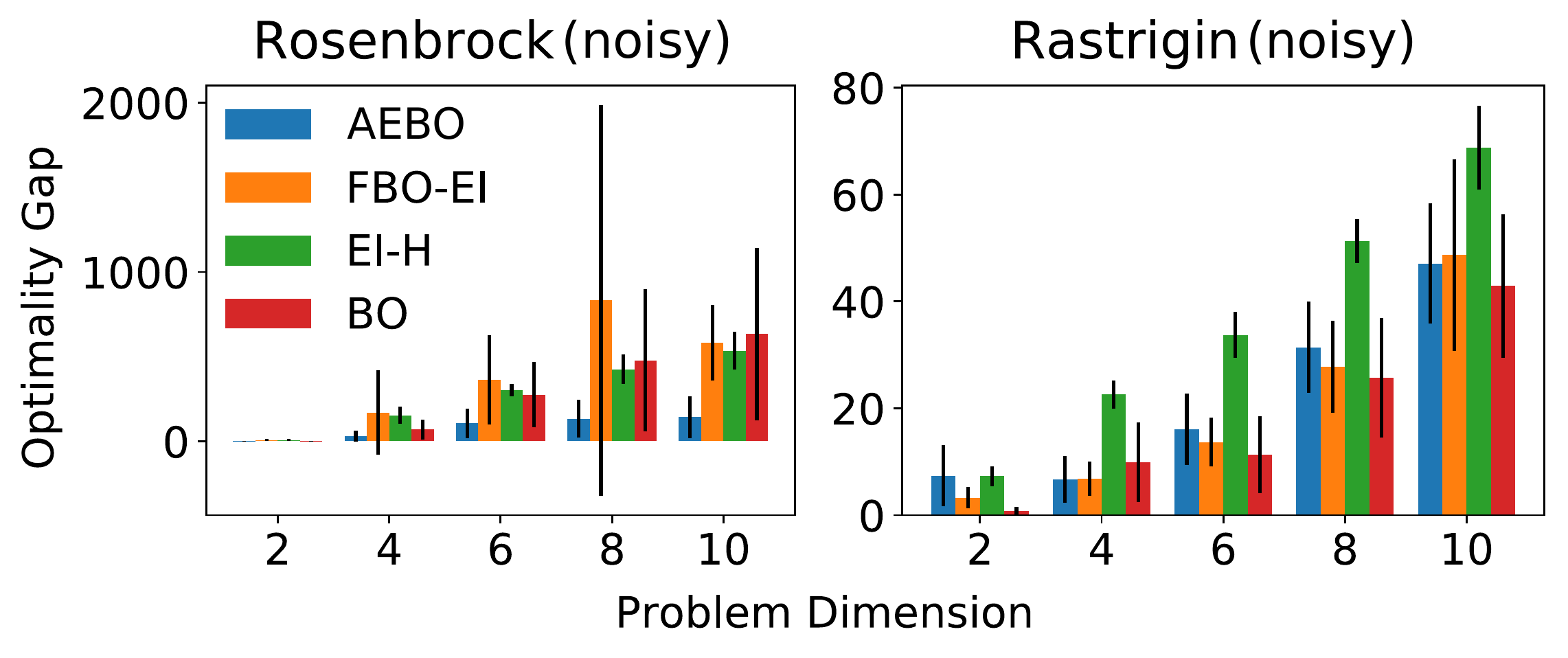}
\caption{Optimality gap results on noisy test functions.}
\label{fig:noise}
\end{figure}

\begin{figure*}[h]
\centering
\includegraphics[width=0.7\textwidth]{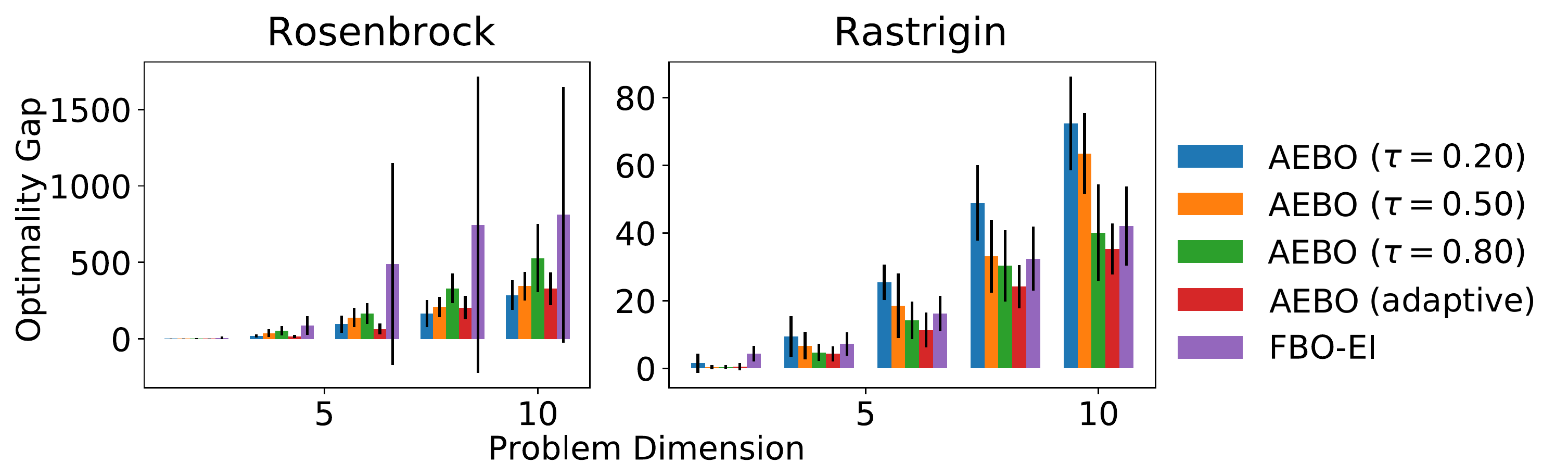}
\caption{The effect of fixed and adaptive $\tau$ on optimality gaps. Here FBO-EI is shown as a baseline.}
\label{fig:tau}
\end{figure*}

We examined the \textit{reachability} of optimal solutions in different methods by measuring the \textit{distance-to-center} metric $s = \|\mathbf{x}^*-\mathbf{c}\|_2$, where $\mathbf{x}^*$ is the optimal solution found and $\mathbf{c}$ is a center point (Fig.~\ref{fig:dim_gap_dist}). In EI-Q and EI-H, $\mathbf{c}$ is the user-specified center; and in AEBO or FBO, $\mathbf{c}$ represents the center of initial bounds. Regularization-based methods like EI-Q and EI-H are biased toward regions near some user-specified center point, thus insufficient exploitation may occur when the optimal solution is far from the center. This is demonstrated in Fig.~\ref{fig:dim_gap_dist}, where EI-H shows a relatively smaller distance-to-center and higher optimality gap comparing to AEBO and FBO.

To evaluate the robustness of AEBO under noise, we tested the case where the observations are corrupted by Gaussian noise with a standard deviation of 0.1. The results are shown in Fig.~\ref{fig:noise}.

To demonstrate the effectiveness of the adaptive exploration-exploitation trade-off, we ran AEBO with both fixed $\tau$ and adaptive $\tau$ solved from Eq.~\ref{eq:set_tau}. The two test functions, Rosenbrock and Rastrigin, have different characteristics and hence prefer different exploration-exploitation trade-offs. Since the Rastrigin function has a large number of local optima, the difficulty for optimizing on Rastrigin is to avoid getting stuck in those local optima. Thus an algorithm with a higher search space expanding rate (\eg, AEBO with a \textit{large} $\tau$) is likely to perform better since it will spend less budget exploiting local optima and more budget expanding towards the global optimum (\ie, \textit{exploration-intense}). In contrast, due to the narrow flat valley in the Rosenbrock function, the difficult part is exploiting near the global optimum to refine the solution. Thus a lower expanding rate (\eg, AEBO with a \textit{small} $\tau$) is likely to be preferred since less budget will be wasted for exploration. The results shown in Fig.~\ref{fig:tau} are consistent with our expectation: the optimality gap increases with the value of the fixed $\tau$ on the Rosenbrock function, while the opposite behavior was observed on the Rastrigin function. However, by using an adaptive $\tau$, AEBO performs better than most other configurations on both test functions. Note that the behavior of FBO-EI is similar to AEBO with a large $\tau$ (without considering the high performance variance on the Rosenbrock function).

As every objective function weights exploitation and exploration differently, BO methods with a fixed expansion schedule may succeed for one function, but fail for another. The AEBO can avoid this by adaptively balancing exploitation-exploration while expanding the search space.

We include the experimental results for constrained BO problems in the supplementary material.

\subsection{MLP on MNIST}

\begin{figure}[h]
\centering
\includegraphics[width=0.5\textwidth]{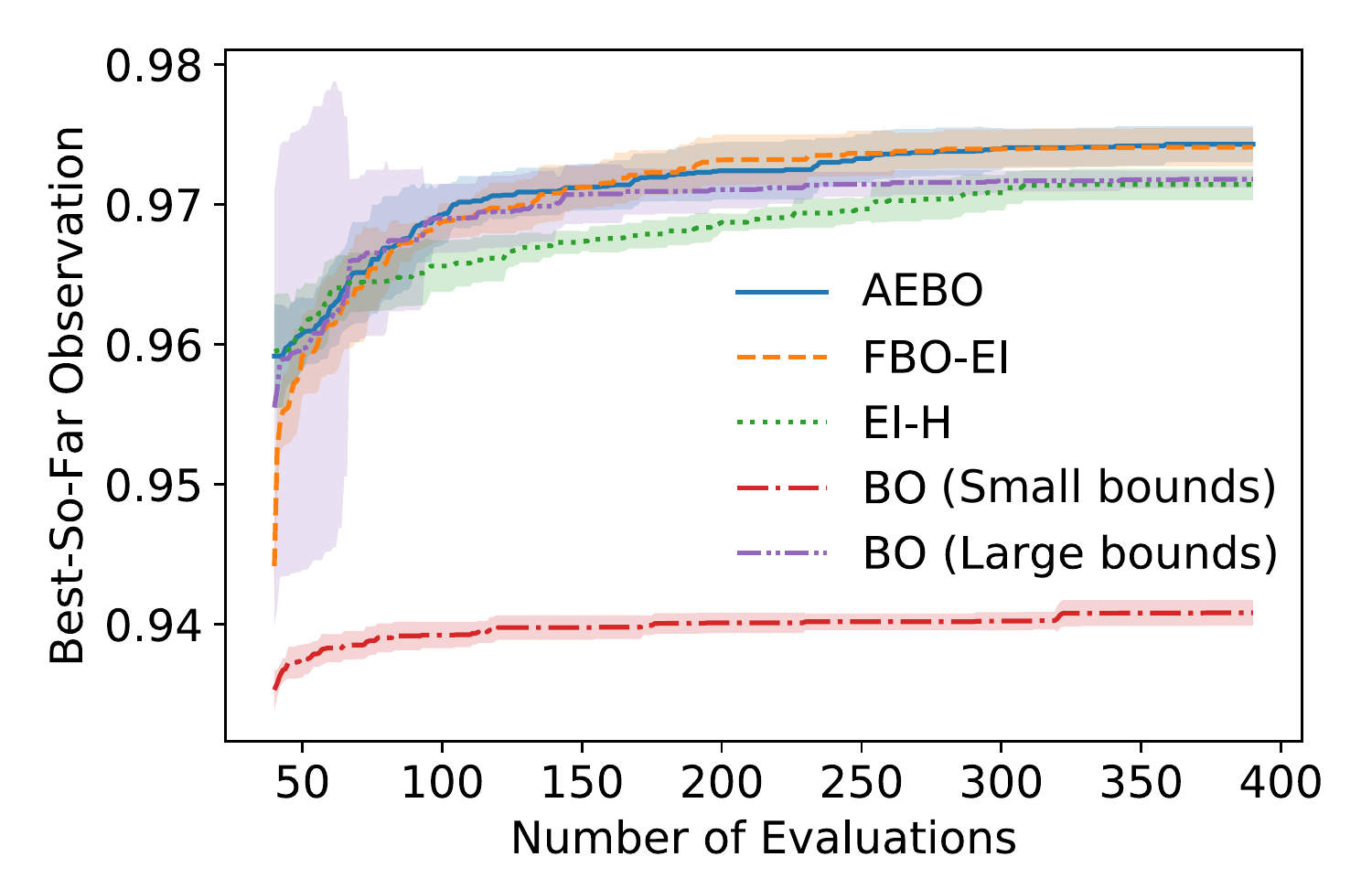}
\caption{Optimization history for the hyperparameter tuning of a MLP trained on MNIST.}
\label{fig:mnist}
\end{figure}

We use the hyperparameter optimization of a multilayer perceptron (MLP) as a real-world example to demonstrate the performance of the AEBO. MNIST was used as the training data. The MLP has three hidden layers, each having 64 hidden units with ReLU activations and was implemented using TensorFlow~\citep{tensorflow2015-whitepaper}. We used Adam~\citep{kingma2014adam} as the MLP's optimizer. We optimized 13 hyperparameters, namely the learning rate, the learning rate decay, the dropout rate for each hidden layer, and the L1 and L2 regularization coefficients for each hidden layer and the output layer. We performed the AEBO, FBO-EI, and EI-H in the log space (base 10) with the initial bounds of $[-5, -4]^7$. We also compared them to standard BO with bounds fixed at $[-5, -4]^7$ (\textit{small bounds}) and $[-10, 1]^7$ (\textit{large bounds}). Note that large and small bounds are different in whether they cover solutions found by unbounded methods, such that the size of the bounds will (for small bounds) or will not (for large bounds) be a factor that limits the solution of BO. The objective is to maximize the accuracy of the MLP. As shown in Fig.~\ref{fig:mnist}, AEBO, FBO-EI, and EI-H found better solutions than BO with small bounds. AEBO out-performed EI-H and BO with large bounds, and was at least as good as FBO-EI.

\section{Conclusion}

We proposed a Bayesian optimization method, AEBO, that gradually expands the search space, so that we can find the global optimum without having to specify the input space bounds that include it. The proposed method only evaluates samples at regions with low GP model uncertainty, and expands the search space adaptively to avoid over-exploration in an expanding search space. This method is useful in cases where we are not confident about the range of the global optimum. The experimental results show that our method out-performs the other state-of-the-art methods in most cases.

In the standard BO, even if the input space bounds are set large enough to cover the global optimum, too much budget may be spent on needlessly exploring the large space. This will result in bad solutions when optimizing an exploitation-intense objective function, as shown by the Rosenbrock and the MLP examples.




\bibliographystyle{spbasic}
\bibliography{refs}

\clearpage
\begin{center}
\textbf{\Large Appendix A: Derivation of the Feasible Domain Bounds}
\end{center}
\setcounter{equation}{0}
\setcounter{figure}{0}
\setcounter{table}{0}
\setcounter{section}{0}
\makeatletter
\renewcommand{\theequation}{A\arabic{equation}}
\renewcommand{\thefigure}{A\arabic{figure}}
\renewcommand{\thetable}{A\arabic{table}}
\renewcommand{\thesection}{A\arabic{section}}
\renewcommand{\citenumfont}[1]{A#1}

In this section, we derive the bounding box that contains the feasible domain of Eq.~\ref{eq:strategy}.

\begin{figure}[h]
\centering
\includegraphics[width=0.25\textwidth]{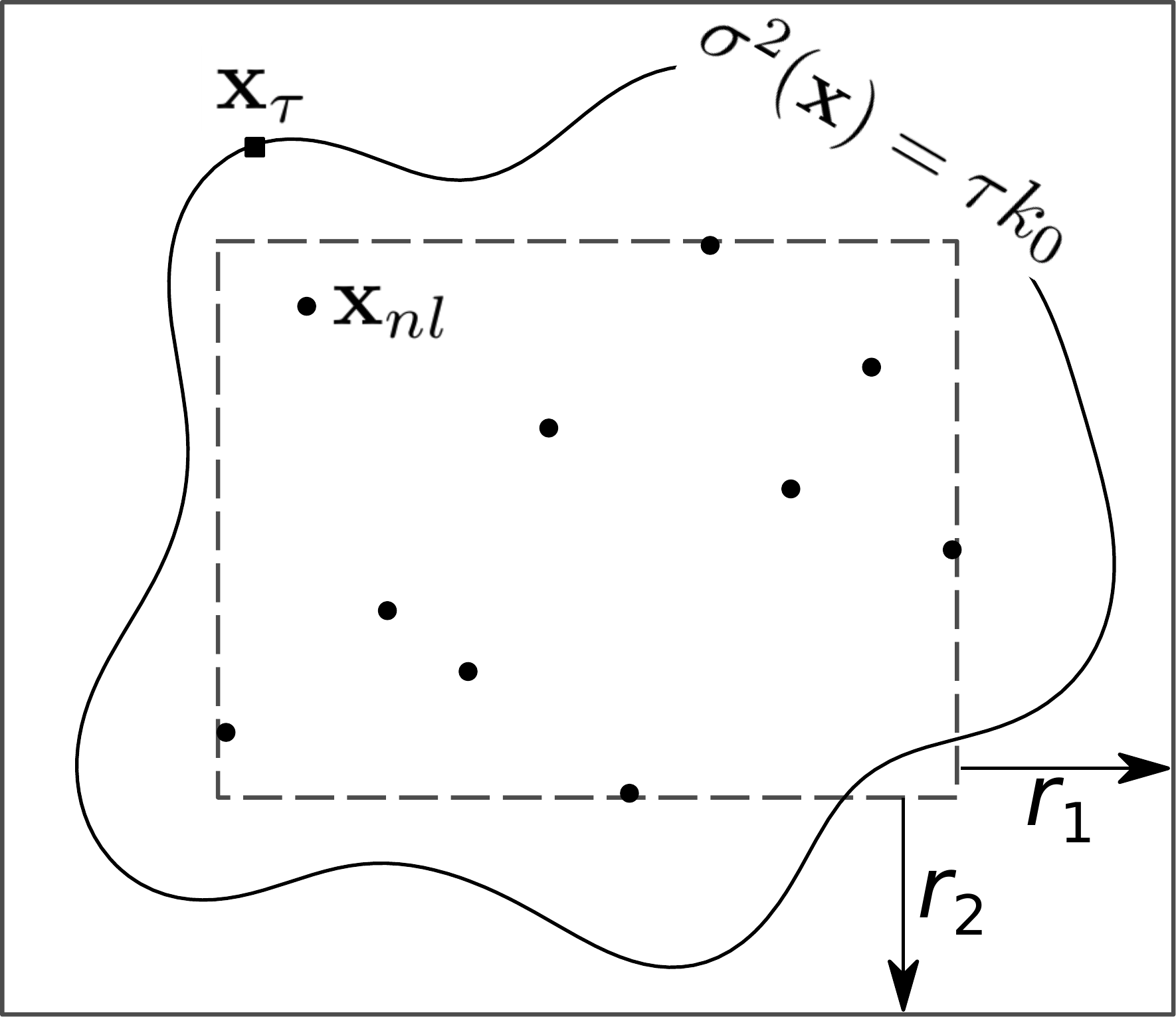}
\caption{Feasible domain bounds. In each iteration, we expand the minimum bounding box of evaluated samples along the $i$-th axis by $r_i$.}
\label{fig:bounds_appendix}
\end{figure}

For any point $\mathbf{x}_\tau$ on the isocontour (Fig.~\ref{fig:bounds_appendix}), \ie, $\sigma^2(\mathbf{x}_\tau) = \tau k_0$, based on Eq.~\ref{eq:variance} we have $k_0-\mathbf{k}_\tau^\top \mathbf{A} \mathbf{k}_\tau = \tau k_0$, or 
\begin{equation}
\mathbf{k}_\tau^\top \mathbf{A} \mathbf{k}_\tau = (1-\tau) k_0
\label{eq:x_tau_0}
\end{equation}
where $\mathbf{A}=(\mathbf{K}+\sigma_n^2\mathbf{I})^{-1}$ and $\mathbf{k}_\tau=\mathbf{k}(\mathbf{x}_\tau)$. 

Since $\mathbf{A}$ is symmetric, we have $\lambda_{\min} \mathbf{k}_\tau^{\top} \mathbf{k}_\tau \leq \mathbf{k}_\tau^\top \mathbf{A} \mathbf{k}_\tau \leq \lambda_{\max} \mathbf{k}_\tau^{\top} \mathbf{k}_\tau$, where $\lambda_{\min}$ and $\lambda_{\max}$ are the smallest and largest eigenvalues of $\mathbf{A}$, respectively. Thus $\mathbf{k}_\tau^{\top} \mathbf{k}_\tau$ have the following bounds for any $\mathbf{x}_\tau$:
\begin{equation}
(1-\tau) k_0 / \lambda_{\max} \leq \mathbf{k}_\tau^{\top} \mathbf{k}_\tau \leq (1-\tau) k_0 / \lambda_{\min}
\label{eq:x_tau_1}
\end{equation}

Suppose $\mathbf{x}_{nl}$ is the nearest evaluated point to $\mathbf{x}_\tau$ (Fig.~\ref{fig:bounds}), the following inequality holds:
\begin{equation}
\mathbf{k}_\tau^{\top} \mathbf{k}_\tau < Nk^2(\mathbf{x}_{nl},\mathbf{x}_\tau)
\label{eq:x_tau_2}
\end{equation}
where $N>1$ is the number of evaluated points.

According to Eq.~\ref{eq:x_tau_1} and Eq.~\ref{eq:x_tau_2}, we have
\begin{equation}
Nk^2(\mathbf{x}_{nl},\mathbf{x}_\tau) > (1-\tau) k_0 / \lambda_{\max}
\label{eq:x_tau_3}
\end{equation}
for any $\mathbf{x}_\tau$. Given any stationary kernel $k(\mathbf{x},\mathbf{x}')=f(\mathbf{\delta})$, where $\mathbf{\delta}=\mathbf{x}-\mathbf{x}'$, we can find the upper bound of $\mathbf{\delta}_i$ for the $i$-th axis. Then we can set that upper bound as the expansion rate $r_i$.

For example, when using the RBF kernel, 
we have
\begin{equation}
k^2(\mathbf{x}_{nl},\mathbf{x}_\tau) = \exp\left(-\frac{1}{2}\sum_{i=1}^d \left(\frac{\delta_i}{l_i}\right)^2\right)
\label{eq:x_tau_4}
\end{equation}
where $\mathbf{\delta} = \mathbf{x}_{nl}-\mathbf{x}_\tau$.
Substituting Eq.~\ref{eq:x_tau_4} into Eq.~\ref{eq:x_tau_3}, we get the following inequality
\begin{equation}
\sum_{i=1}^d \frac{\delta_i^2}{Cl_i^2} < 1
\label{eq:x_tau_5}
\end{equation}
where $C=-\log((1-\tau)k_0/(N\lambda_{\max}))$. Equation~\ref{eq:x_tau_5} shows that $\mathbf{x}_\tau$ is inside a $d$-dimensional hyperellipsoid that centered at $\mathbf{x}_{nl}$ with $r_i=\sqrt{C}l_i$ corresponding to half the length of the $i$-th principal axis. 

By setting the bounds of the $i$-th dimension as $\left[\min_j\{x_i^{(j)}\}-r_i, \max_j\{x_i^{(j)}\}+r_i\right]$, we can include the entire feasible domain of Eq.~\ref{eq:strategy}. This means that in each iteration, we get the minimum bounding box of all evaluated samples, and expand the bounding box along the $i$-th axis by $r_i$ (Fig.~\ref{fig:bounds}). Then constrained global optimization of the acquisition function is performed within the new bounds.

In practice, because Eq.~\ref{eq:x_tau_2} is usually quite loose (especially when $N$ is large), the above derived bounds are usually unnecessarily large, causing large volume of infeasible domain inside the bounds. In that case, we can replace $\lambda_{\max}$ with $\lambda_{\min}$ in Eq.~\ref{eq:x_tau_3}, \ie, substituting the upper bound of $\mathbf{k}_\tau^{\top} \mathbf{k}_\tau$ (Eq.~\ref{eq:x_tau_1}) into Eq.~\ref{eq:x_tau_3}.

\begin{center}
\textbf{\Large Appendix B: Experiments for Constrained BO Problems}
\end{center}
\setcounter{equation}{0}
\setcounter{figure}{0}
\setcounter{table}{0}
\setcounter{section}{0}
\makeatletter
\renewcommand{\theequation}{B\arabic{equation}}
\renewcommand{\thefigure}{B\arabic{figure}}
\renewcommand{\thetable}{B\arabic{table}}
\renewcommand{\thesection}{B\arabic{section}}
\renewcommand{\citenumfont}[1]{B#1}

\begin{figure*}[h]
\centering
\includegraphics[width=0.9\textwidth]{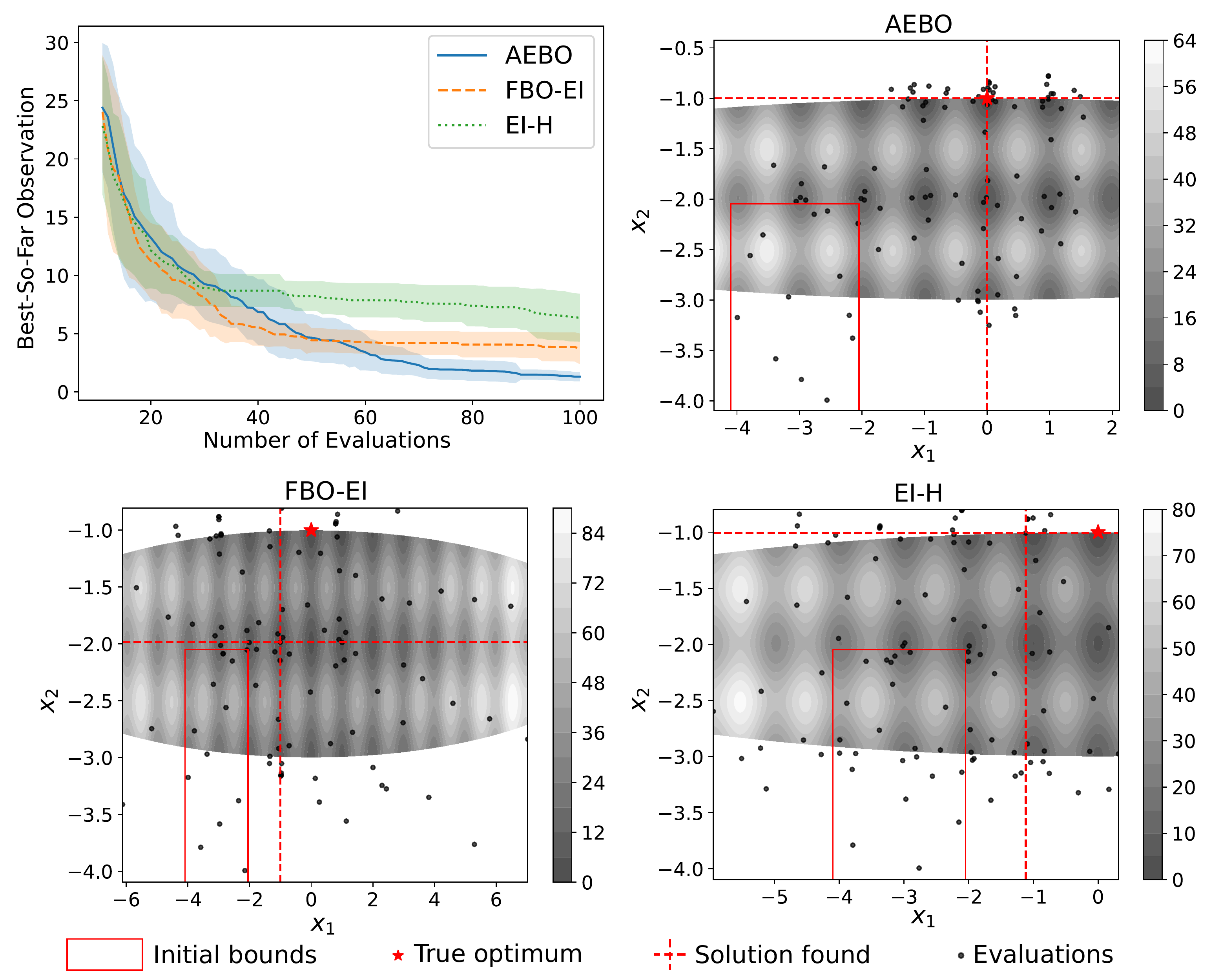}
\caption{Optimization history and evaluated points for the constrained Rastrigin problem}
\label{fig:constrained_rastrigin}
\end{figure*}

We created two test problems to evaluate the performance of AEBO in dealing with constrained BO problems. Specifically, the constrained Rastrigin problem uses the Rastrigin function as the objective function, and the feasible domain is defined by an ellipse $0.01x_1^2 + (x_2+2)^2 \leq 1$ (Fig.~\ref{fig:constrained_rastrigin}). The Nowacki beam problem is a real-world test problem originally described by Nowacki~\citep{nowacki1980modelling,singh2017sequential}. The goal is to minimize the cross-sectional area of a tip-loaded cantilever beam subject to certain constraints.\footnote{The original problem is a multi-objective optimization problem that minimizes both the cross-sectional area and the bending stress. Here we only consider the first objective and limit the second objective (\ie, the bending stress should be smaller than the yield stress of the material) to form an extra constraint.}

\begin{figure*}[h]
\centering
\includegraphics[width=0.9\textwidth]{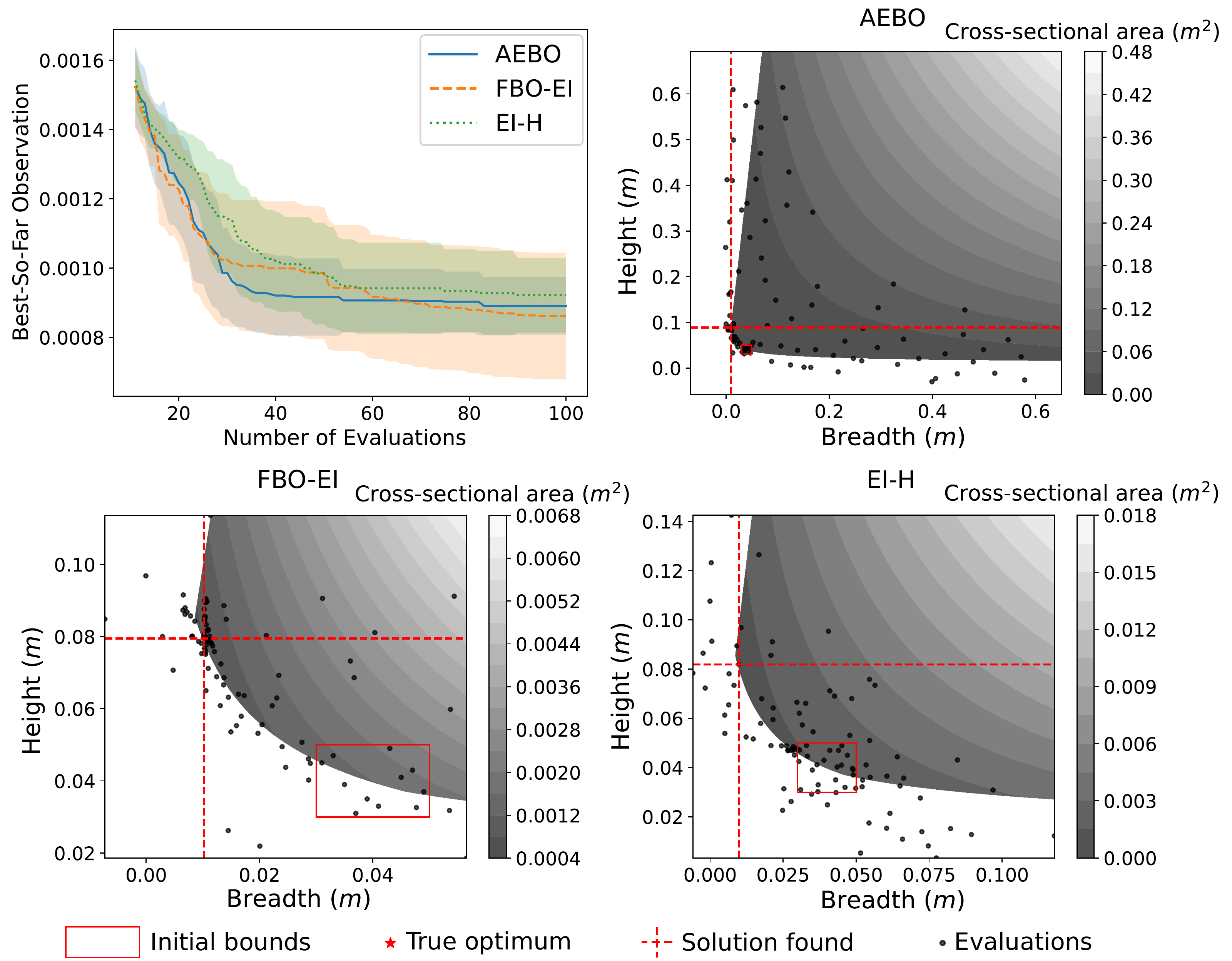}
\caption{Optimization history and evaluated points for the Nowacki beam problem}
\label{fig:beam}
\end{figure*}

The results of the two problems are shown in Figs.~\ref{fig:constrained_rastrigin} and \ref{fig:beam}. For the constrained Rastrigin problem, AEBO achieved a better solution than the other two methods. For the Nowacki beam problem, FBO-EI's solution has the lowest mean value but a very high variance; while AEBO found a fairly close optimal solution with much lower variance. The evaluated points by AEBO were dense near optima (either local or global). This behavior was, however, not obvious for the other two methods. This is likely because that FBO-EI and EI-H over-trusted the GP posterior even where its uncertainty was high. This resulted in sampling patterns with too much randomness, and hence higher variance of optimal solutions.

\end{document}